\documentclass[mathematics,article,accept,moreauthors,pdftex]{Definitions/mdpi}

\firstpage{1} 
\pubvolume{1}
\issuenum{1}
\articlenumber{0}
\pubyear{2021}
\copyrightyear{2020}
\datereceived{} 
\dateaccepted{} 
\datepublished{} 
\hreflink{https://doi.org/} % If needed use \linebreak
\usepackage{amsfonts}
\usepackage{amssymb}
\usepackage[utf8]{inputenc}
\usepackage[english]{babel}
\usepackage{array}
\usepackage{subcaption}
\usepackage{lmodern}

\def\EI{\mathcal{E}_\infty}
\def\R{\mathbb R}

\def\e{\varepsilon}
\def\t{\mathfrak{t}}
\def\tT{(\t,T)}
\def\r{\mathfrak{r}}
\def\f{\mathfrak{f}}
\def\g{\mathfrak{g}}
\def\cl{\overline}
\def\M{\mathcal{M}}
\def\m{{\scriptstyle{\mathpzc{M}}}}

\newcommand\sbullet[1][.5]{\mathbin{\vcenter{\hbox{\scalebox{#1}{$\bullet$}}}}}
\DeclareMathAlphabet{\mathpzc}{OT1}{pzc}{m}{it}

\Title{TAC Method for Fitting Exponential Autoregressive Models and Others: Applications in Economy and Finance}

\TitleCitation{TAC and Applications in Economy and Finance}

\Author{Javier Cabello Sánchez $^{1,\dagger,\ddagger}$\orcidA{}*, Juan Antonio Fernández Torvisco $^{2,\ddagger}\orcidB{}$ and Mariano R. Arias $^{3,\ddagger}$\orcidC{}*}

\AuthorNames{Javier Cabello~Sánchez, Juan Antonio Fernández~Torvisco and Mariano R. Arias}

\AuthorCitation{Cabello~Sánchez, J.; F.~Torvisco, J.A.; R.~Arias, M.}
\address[1]{%
$^{1}$ \quad Facultad de Ciencias, Universidad de Extremadura, Avda. de Elvas s/n, 06006 Badajoz. Spain; coco@unex.es\\
$^{2}$ \quad Facultad de Ciencias, Universidad de Extremadura, Avda. de Elvas s/n, 06006 Badajoz. Spain; jfernandck@alumnos.unex.es\\
$^{3}$ \quad Facultad de Ciencias, Universidad de Extremadura, Avda. de Elvas s/n, 06006 Badajoz. Spain; arias@unex.es}

\corres{Correspondence: coco@unex.es; arias@unex.es}

\firstnote{Current address: Facultad de Ciencias, Universidad de Extremadura, Avda. de Elvas s/n, 06006 Badajoz. Spain.} 
\secondnote{These authors contributed equally to this work.}
\abstract{There are a couple of purposes in this paper: to study a problem of approximation with exponential functions and to show its relevance for the economic science.
We present results that completely solve the problem of the best approximation by means of exponential functions and we will be able to determine what kind of data is suitable to be fitted. Data will be approximated using TAC (implemented in the R-package \textit{nlstac}), a numerical algorithm for fitting data by exponential patterns without initial guess designed by the authors. We check one more time the robustness of this algorithm by successfully applying it to two very distant areas of economy: demand curves and nonlinear time series. This shows TAC's utility and highlights how far this algorithm could be used.}


\begin{document}
%%%%%%%%%%%%%%%%%%%%%%%%%%%%%%%%%%%%%%%%%%
%\setcounter{section}{-1} %% Remove this when starting to work on the template.

\section{Introduction}
This paper is planned to cover a couple of major objectives. Broadly, the first 
one is to solve one of the remaining issues in~\cite{TAC1}. Later, in this same 
introduction, this first objective will be introduced in a detailed way. The 
second one is to exemplify the interest, for the economic science, of the algorithm 
presented in~\cite{TAC1} and implemented in the R-package \textit{nlstac}, 
see~\cite{nlstac}. We have developed this algorithm in order to fit data 
coming from an exponential decay. Therefore we will illustrate the interest 
of this algorithm by fitting the pattern in a couple of cases related with 
different economic problems. 

The first economic problem deals with demand curves. Economic demand curves can 
be used to map the relationship between the consumption of a good and its price. 
When plotting price (actually the logarithm of the price) and consumption, we 
obtain a curve with a negative slope, meaning that when price increases, demand 
decreases. Hursh and Silberbeg proposed in~\cite{example_hursh_economic_demand} 
an equation to model this situation; in this paper, we will fit some data using 
that model.

The second economic problem is about nonlinear time series models. Many financial 
time series display typical nonlinear characteristics, so many authors such as 
\cite{example_Hans_non_linear} consider to apply nonlinear models. Although our 
algorithm was not designed for this kind of problems, we will obtain good results 
using it. In this example we will focus on the model that, among all nonlinear 
time series models, seems to have more relevance in the literature, namely, the 
exponential autoregressive model.

Before we get into the first objective, let us see the structure of this paper. 
Section \ref{s_small} deals with approximations by means of exponential functions 
measuring the error with the $\max$-norm when fitting a \textit{small} set of data: 
3 or 4 observations. Later, Section \ref{s_symmetric} is devoted to some symmetric 
cases that could happen. Section \ref{s_general} deals with approximation in 
\textit{general} datasets. 
Section~\ref{s_examples} gathers two examples about 
Newton Law of Cooling and directly apply what have been developed in this paper. 
Section~\ref{s_examples_economy} shows examples related to economy and use the 
R-package \textit{nlstac} for the calculations. Although the economical section represents 1 out of the 6 sections of this paper, it remains a very important one, and everything we have been working on before directly apply there.

\begin{figure}
\begin{center}
\includegraphics[height=0.2\paperheight]{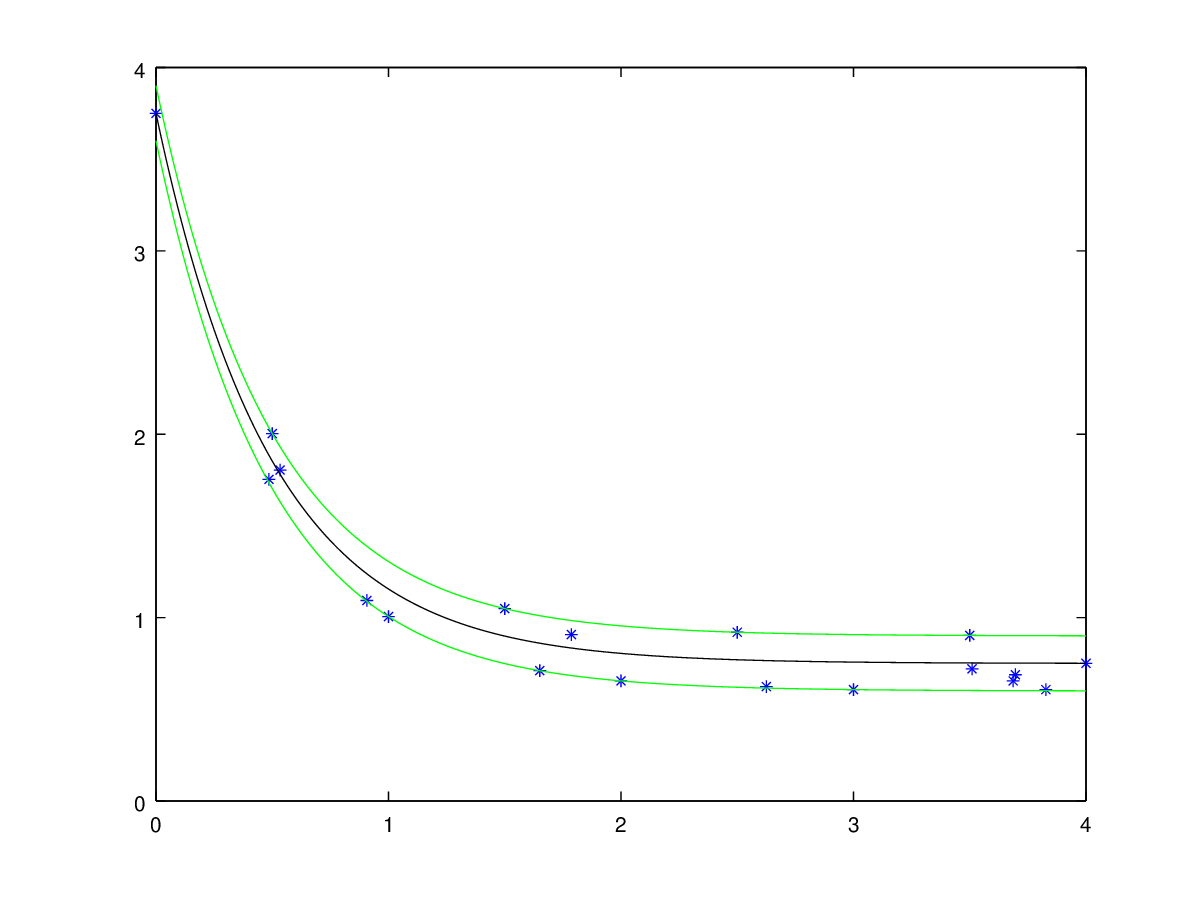}
\end{center}
\caption{In blue, the points, in black the best approximation and in green the 
upper and lower borders of the narrowest band that contains $(\t,T)$. 
Yes, the band has constant width. }\label{DibujoBanda}
\end{figure}

We will focus in the first objective in Sections \ref{s_small} to \ref{s_examples}. 
In those sections we will deal with approximations by means of exponential 
functions and we will measure the error with the $\max$-norm. 
So, when we say that some function $f$ is the best approximation for some 
data $(\t,T)$, we will mean that we have 
$\t=(t_1\ldots,t_n), T=(T_1\ldots,T_n)\in\R^n$, 
with $t_i<t_{i+1}$ for every $i$, and that $\{(t,f(t)):t\in[t_1,t_n]\}\subset \R^2$ 
is the centre of the narrowest band that contains every point and has exponential 
shape, see Figure~\ref{DibujoBanda}. 

Before we go any further, let us introduce the definition of quasiconvex function. 
A real function defined in a linear topological space $X$ is quasiconvex whenever it fulfills 
\begin{equation*}
f(\lambda x + (1-\lambda) y) \leq \max\{f(x), f(y)\}, \forall x,y \in X 
\text{ and } \lambda \in (0,1).
\end{equation*} 
This definition can be consulted in, for example, 
\cite{TAC1} or \cite{harvey1971review}.

The authors already proved in~\cite{TAC1} that for every $k\in\R,\ k\neq 0$, there 
exists the best approximation $f_k(t)=a_k\exp(kt)+b_k$ amongst all the 
functions of the form $a\exp(kt+b)$, with $a, b\in\R$. We also showed that, 
given any dataset $(\t,T)$, the function $\EI(k)$ that assigns to every 
$k\in(-\infty,0)$ the error $\max|T_i-f_k(t_i)|$ is quasiconvex, so there are 
two options. Namely, either $\EI$ attains its minimum at some $k$ or it is monotonic. 
If $\EI$ is monotonic, then either it is increasing and the minimum would 
be attained, so to say, at $-\infty$, or it is decreasing and attains its 
minimum at $k=0$. We will also study what happens for positive $k$, so 
we will need to pay attention not only to the behaviour of the exponentials but 
also to their limits when $k\to 0$ and $k\to\pm\infty$, see 
Proposition~\ref{limites} and Subsection~\ref{sbs_limit}. 

\begin{Notation}\label{Desglose}
Our main results show that every dataset $(\t,T)$, with at least four data, 
fulfils one of the following conditions:
\begin{itemize}
\item There exists one triple $(a, b, k)\in\R^3$ with $ak\neq 0$ such that 
$a\exp(kt)+b$ is the best possible approximation, i.e., 
\begin{equation*}
\max_{i=1,\ldots,n}\{|a\exp(kt_i)+b-T_i|\}<\max_{i=1,\ldots,n}\{|f(t_i)-T_i|\}
\end{equation*}
whenever $f:\R\to\R$ is a different exponential. 
\item There are two indices $i_1<i_3$ where $\max\{T_1,\ldots,T_n\}$ is attained 
and there is some $i_2$, with $i_1<i_2<i_3$, such that $\min\{T_1,\ldots,T_n\}=T_{i_2}$. In this 
case, the best approximation by means of exponentials does not exist and the constant 
$\frac 12(T_{i_1}+T_{i_2})$ approximates $(\t,T)$ better than any strictly 
monotonic function --in particular any exponential. So, for every $k\in\R$, 
the best approximation with the form $a\exp(kt)+b$ has $a=0$ and 
$b=\frac 12(T_{i_1}+T_{i_2})$. It happens exactly the same when the maximum 
is attained at $i_2$ and the minimum at $i_1$ and $i_3$, with $i_1<i_2<i_3$. 
In both cases, the function $\EI$ is constant. 
\item $T_1=\max\{T_1,\ldots,T_n\}$, $T_{i_2}=\min\{T_1,\ldots,T_n\}$ and 
$T$ attains its second greatest value at $i_3>i_2$. The best approximation by means 
of exponentials does not exist and every exponential approximates $(\t,T)$ worse 
than any function fulfiling $f(t_1)=T_1+\frac 12(T_{i_2}-T_{i_3}), 
f(t_i)=\frac 12(T_{i_2}+T_{i_3})$ for $i\geq 2$. The pointwise limit of the 
best approximations when $k\to-\infty$ takes these values. (The symmetric cases 
belong to this kind of limits, maybe with $k\to\infty$ instead of $k\to-\infty$). 
If this happens, $\EI$ is increasing in $(-\infty,0).$
\item There are some $c,d\in\R$ such that the line $ct+d$ approximates $(\t,T)$ better 
than any exponential. In this case, each $ct_i+d$ is the limit when $k\to 0$ of 
the values in $t_i$ of the best approximations with $k$ as exponent. This happens  
when there are four indices $i_1<i_2<i_3<i_4$ or $i_2<i_1<i_4<i_3$ such that 
\begin{equation*}
T_{i_j}-(ct_{i_j}+d)=(-1)^j\max |T_{i}-(ct_{i}+d)|, j= 1, 2, 3, 4.
\end{equation*}
This implies that $\EI$ decreases in $(-\infty,0).$
\end{itemize}
\end{Notation}

\begin{Remark}
What happens with this kind of functions is the following: 
consider two exponentials that agree at $\alpha<\beta\in\R$, 
say $g(\alpha)=f(\alpha)>g(\beta)=f(\beta)$. Then, the following are equivalent: 
\begin{itemize}
\item $g(\gamma)<f(\gamma)$ for some $\gamma\in(\alpha,\beta)$. 
\item $g(\gamma)<f(\gamma)$ for every $\gamma\in(\alpha,\beta)$. 
\item $g(\gamma)>f(\gamma)$ for some $\gamma\in(-\infty,\alpha)\cup(\beta,\infty)$. 
\item $g(\gamma)>f(\gamma)$ for every $\gamma\in(-\infty,\alpha)\cup(\beta,\infty)$. 
\end{itemize}
A visual way to look at this is the following. 
Consider a wooden slat supported on two points and imagine we put a load
between the supports. When we increase the load, the slat lowers between the 
two supports but the other part of the slat raises. 
For these functions, the behaviour is similar --if two of them agree at 
$\alpha$ and $\beta$ then one function is greater than 
the other in $(\alpha,\beta)$ and lower outside $[\alpha,\beta]$. 

Besides, if $\alpha<\beta$ and $f(\beta)<f(\alpha)$ then for each $(\gamma,x)$ 
that does not lie in the line defined by $(\alpha,f(\alpha))$ and 
$(\beta,f(\beta))$ and belongs to 
\begin{equation}\label{set}
\big((-\infty,\alpha)\times(f(\alpha),\infty)\big)
\cup\big((\alpha,\beta)\times(f(\beta),f(\alpha))\big)
\cup\big((\beta,\infty)\times(-\infty,f(\beta))\big)
\end{equation}
there is exactly one exponential $h$ such that $h(\alpha)=f(\alpha)$, 
$h(\beta)=f(\beta)$ and $h(x)=\gamma$. Of course, if $(\gamma,x)$ does not 
belong to the set given by~(\ref{set}), then there is no monotonic function that fulfils 
the later. The existence of such an exponential is a straightforward 
consequence of~\cite[Lemma 2.10]{TAC1} --we will develop this later, see 
Proposition~\ref{decinchom}. 
\end{Remark}

\begin{Remark}
A significant trouble when dealing with the problem of approximating datasets 
with exponentials has been to find conditions determining whether some dataset 
is worth the try or not. The only way we have found to answer this problem has 
been to identify the most general conditions that ensure that some dataset has 
{\em one} best approximation by exponentials --needless to say, this has been 
a very sinewy problem. 
The different behaviours described in Notation~\ref{Desglose} can give a hint about 
the several different details that we will need to deal with, but there is still 
some casuistry that we need to break down. Namely, our main interest in these 
results comes from the fact that they can be applied to exponential decays, 
which appear in several real-life problems --the introduction in~\cite{TAC1} 
presents quite a few examples. The typical data that we have worked 
with are easily recognizable, but we needed to determine when the data may be fitted with a 
decreasing, convex, function --like the exponentials $f(t)=a\exp(kt)+b$ with 
$a>0, k<0$. The easiest way we have found is:

\begin{itemize}
\item If some data $(\t,T)$ are to be fitted with a decreasing function and we are 
measuring the error with the $\max$-norm, then the maximum value in $T$ must be attained 
before the minimum. There may exist more than one index where they are attained, but 
every appearance of the maximum must lie before every appearance of the minimum. 
In short, if $T_i=\max(T)$ and $T_j=\min(T)$, then $i<j$. 
\item Moreover, if we are going to approximate $(\t,T)$ with a convex function, 
the dataset must have some kind of {\em convexity}. 
The only way we have found to state this is: \\
$\heartsuit$
``Let $r(t)$ be the line that best approximates $(\t,T)$. Then 
$(T_1-r(t_1),\ldots,T_n-r(t_n))$ has two maxima and one minimum between them." \\
Thanks to Chebyshev's Alternation Theorem (the polynomial $p_n$ is the best approximation of function $f$ in $[a,b]$ if and only if there exists $n+2$ points $a\leq t_1<t_2< \ldots < t_{n+2} \leq b$ where $\vert f(t)-p_n(t) \vert$ attains its maximum and $f(t_i)-p_n(t_i)=p_n(t_{i+1})-f(t_{i+1})$, see for example \cite[Theorem 8 page 29]{lorentz_approximation} or \cite{poreda_chebyshev}) we know that the line that best approximates any dataset behaves either this way, 
either the opposite or the third case in~\ref{Desglose}. Please observe that this 
theorem would not apply so easily to approximations with general degree polynomials. 
\end{itemize}
\end{Remark}

%\subsection{Notations}

Before we go any further, let us comment something about the notation that will be used. For the remaining of the paper we will always take $n$ as the number of 
coordinates of $\t$ and $T$, i.e., $\t, T\in\R^n$. Besides, $\t\in\R^n$ will fulfil $t_1<t_2<\ldots<t_n$. 

Moreover, for any $k\neq 0$ we will always denote as $f_k(t)=a_k\exp(kt)+b_k$ 
the best approximation with the form $f(t)=a\exp(kt)+b$ with $a,b\in\R$. 
%\end{Notation}

%\begin{Notation}\label{fraktur}
In a try to ease the notation, whenever we have some function $f:\R\to\R$ 
and $\t=(t_1,\ldots,t_n)\in\R^n$, $\f(\t)$ will denote 
$(f(t_1),\ldots,f(t_n))\in\R^n$, and the same will apply to any 
{\em fraktur} character: $\g(\t), \r(\t), \f_1(\t), \f_2(\t), \f_k(\t)$ would represent 
the same for $g, r, f_1, f_2, f_k:\R\to\R$ and so on. 

Given any vector $v=(v_1,\ldots,v_n)$, $\M(v)$ will denote its maximum 
and $\m(v)$ will denote its minimum. 
%\end{Notation}

\section{Small datasets}\label{s_small}
In this Section we are going to show some not too complicated, general results 
about the behaviour of exponentials that will allow us to prove our main results
in Section~\ref{s_general}. We will focus only in the approximation of the most 
simple datasets, with $n=3$ or $n=4$. 

We will begin with Proposition~\ref{rrr}, just a slight modification 
of~\cite[Proposition 2.3]{TAC1} that will be useful for the subsequent results.
Later, in Lemma~\ref{fijok3r}, we will find the expression of the best approximation for 
$((t_1,t_2,t_3),(T_1,T_2,T_3))$ for each fixed $k$ (please note that with $n=2$, 
for every $k$ there is an exponential that interpolates the data). 
In Lemma~\ref{44} we study the case $n=4$, determining a technical condition on 
$(\t,T)$ that ensures that the best approximation exists and it is unique and, 
moreover, we kind of determine analytically this best approximation. 
 
\begin{Proposition}[\cite{TAC1}]\label{rrr}
Let $k\neq 0$, $a_k, b_k\in\R$ such that $\f_k(\t)=a_k\exp(k\t)+b_k$ is the best 
approximation to $T$ for this $k$, i.e.,
\begin{equation}\label{eqrrr}
\|\f_k(\t)-T\|_\infty=\min\{\|a\exp(k\t)+b-T\|_\infty:a, b\in\R\}. 
\end{equation}
Then, there exist indices $1\leq i<j<m\leq n$ such that 
$f_k(t_i)-T_i=T_j-f_k(t_j)=f_k(t_m)-T_m=\pm \|\f_k(\t)-T\|_\infty$. 

Reciprocally, if $a_k$ and $b_k$ fulfil this condition, then $a_k\exp(k\t)+b_k$ 
is the best approximation to $T$ for this $k$. 
\end{Proposition}

\begin{Lemma}\label{fijok3r}
Let $n=3$, $k\neq 0$ and $T_1\neq T_3$. Then, the best approximation to 
$\tT$ by means of exponentials has these coefficients:  
\begin{equation}\label{ak3bk3}
a_k=\frac{T_1-T_3}{\exp(kt_1)-\exp(kt_3)}, b_k=\frac 12(T_1-a_k\exp(kt_1)+T_2-a_k\exp(kt_2)). 
\end{equation}
\end{Lemma}

\begin{proof}
It is clear that $T_1-f(t_1)=T_3-f(t_3)$, and a simple computation shows that 
$T_1-f(t_1)=f(t_2)-T_2$ also holds. Indeed, 
\begin{equation*}
\begin{split}
f(t_1)& -f(t_3)=a_k\exp(kt_1)+b_k-a_k\exp(kt_3)-b_k=%\\ &
a_k(\exp(kt_1)-\exp(kt_3))=\\ &
\frac{T_1-T_3}{\exp(kt_1)-\exp(kt_3)}(\exp(kt_1)-\exp(kt_3))=T_1-T_3. \\
f(t_1)& +f(t_2)=a_k\exp(kt_1)+b_k+a_k\exp(kt_2)+b_k=
\\ & a_k(\exp(kt_1)+\exp(kt_2))+ 2\frac 12(T_1-a_k\exp(kt_1)+T_2-a_k\exp(kt_2))
=T_1+T_2. 
\end{split}
\end{equation*}
By Proposition~\ref{rrr}, this is enough to ensure that $a_k$ and $b_k$ are optimal. 
\end{proof}

\begin{Remark}
Please observe that $a_k$ does not depend on $(t_2,T_2)$. 
\end{Remark}

\begin{Lemma}\label{44}
Let $n=4$ and $\tT$ such that 
$(T_1-T_3)/(t_3-t_1)>(T_2-T_4)/(t_4-t_2)>0$. Then, there exists a 
unique exponential $f(t)=a\exp(kt)+b$, with $k<0, a>0$ and $b\in\R$ such that 
\begin{equation}\label{Eq44}
T_1-f(t_1)=-(T_2-f(t_2))=T_3-f(t_3)=-(T_4-f(t_4)).
\end{equation}
Moreover, this exponential is the best approximation to $(\t,T)$. 
\end{Lemma}

\begin{proof}
Let $\tT$ be as in the statement. 
For $k\in\R, k\neq 0$, there exist unique $a, b\in\R$ such that 
$a\exp(kt_1)+b=T_1$ and $a\exp(kt_3)+b=T_3$. Namely, $a$ is as 
in~\eqref{ak3bk3} and $b=T_1-a\exp(kt_1)$. 

Indeed, $T_1-f(t_1)=T_3-f(t_3)$ means that $T_1-(a\exp(kt_1)+b)=T_3-(a\exp(kt_3)+b)$, so 
$T_1-T_3=a(\exp(kt_1)-\exp(kt_3))$ and each $k\neq 0$ determines 
$a^+_k=\frac{T_1-T_3}{\exp(kt_1)-\exp(kt_3)}$. 
The same way, $T_2-f(t_2)=T_4-f(t_4)$ determines $a^-_k=\frac{T_2-T_4}{\exp(kt_2)-\exp(kt_4)}$. 

So, the equalities \eqref{Eq44} hold if and only if, for some $k\neq 0$, we have 
$a^+_k=a^-_k$. Equivalently, 
\begin{equation}
\frac{\exp(kt_1)-\exp(kt_3)}{T_1-T_3}=\frac{\exp(kt_2)-\exp(kt_4)}{T_2-T_4}. 
\end{equation}
Please observe that this equality holds trivially when $k=0$ and that, as both 
$T_1-T_3$ and $T_2-T_4$ are positive, we are not trying to divide by 0. 

If we put $z=\exp(k)$, the last equality can be written as 
\begin{equation*}
(T_1-T_3)z^{t_4}-(T_2-T_4)z^{t_3}-(T_1-T_3)z^{t_2}+(T_2-T_4)z^{t_1}=0.
\end{equation*}
We will denote as $p(z)$ the left hand side of this equality. 

As we are only interested in positive roots of $p$, we can divide by $z^{t_1}$ 
and consider $q(z)=(T_1-T_3)z^{s_4}-(T_2-T_4)z^{s_3}-(T_1-T_3)z^{s_2}+T_2-T_4,$
with $s_i=t_i-t_1$ for $i=2, 3, 4$. 

Taking into account that $q(0)=T_2-T_4>0$, that $q$ obviously vanishes at 1 (but this 
root corresponds to the void case $k=0$ and so $a^+_k$ and $a^-_k$ are not defined) 
and also that the limit of $q(z)$ is $\infty$ as $z$ goes to $\infty$, 
there must exist another $z\in(0,\infty)$, maybe $z=1$, such that $p(z)=q(z)=0$. 
By Descartes' rule of signs, see~\cite{Haukkanen} Theorem~2.2, both $q$ and $p$ 
have at most two positive roots, so there is exactly another positive root of $p$. 
To determine whether this root is greater or smaller than 1, 
we can compute the derivative of $p$ at 1. 
\begin{equation*}
p'(z)=t_4(T_1-T_3)z^{t_4-1}-t_3(T_2-T_4)z^{t_3-1}-t_2(T_1-T_3)z^{t_2-1}+t_1(T_2-T_4)z^{t_1-1},
\end{equation*}
so $p'(1)=(t_4-t_2)(T_1-T_3)-(T_2-T_4)(t_3-t_1)$. This is positive provided  
$(T_1-T_3)/(t_3-t_1)>(T_2-T_4)/(t_4-t_2)>0$, so the other root of $p$ lies between 
0 and 1 whenever the condition in the statement is fulfiled. 

So, there exists just one $k=\log(z)\in(-\infty,0)$ for which 
\begin{equation*}
\frac{T_1-T_3}{\exp(kt_1)-\exp(kt_3)}=\frac{T_2-T_4}{\exp(kt_2)-\exp(kt_4)}. 
\end{equation*}
Now, taking 
\begin{equation*}
a=\frac{T_1-T_3}{\exp(kt_1)-\exp(kt_3)}=\frac{T_2-T_4}{\exp(kt_2)-\exp(kt_4)}, 
b=\frac 12(T_1-a\exp(kt_1)+T_2-a\exp(kt_2))
\end{equation*}  
and $f(t)=a\exp(kt)+b$ we have the function we were looking for. 

As for the moreover part, suppose that there exist $\cl{a}, \cl{b}, \cl{k}$ such 
that $\cl{f}(t)=\cl{a}\exp(\cl{k}t)+\cl{b}$ approximates $(\t,T)$ at least as 
well as $f$. We may suppose that $T_1-f(t_1)=f(t_2)-T_2=T_3-f(t_3)=f(t_4)-T_4=r>0$. 
Now, the conditions for $\cl{f}$ can be rewritten as 
\begin{equation*}
\cl{f}(t_1)\geq f(t_1), \cl{f}(t_2)\leq f(t_2), \cl{f}(t_3)\geq f(t_3), \cl{f}(t_4)\leq f(t_4).
\end{equation*}
By~\cite[Lemma 2.8]{TAC1}, this means that $\cl{f}=f$. 
\footnote{¿No sería más directo decir que por el corolario 2.9 de TAC1?}
\end{proof}

\begin{Remark}
This Lemma gives a kind of analytic solution to the best approximation problem, 
with the only obstruction of being able to determine the other root of $p$. 
In the next section, we do the same with the symmetric cases and give actual 
analytic solutions to the same problem when the data are not good to be 
approximated by exponentials, ironically. 
\end{Remark}

\section{Symmetric cases and limits}\label{s_symmetric}
In this section we will focus in those cases that do not match with the problem 
we have in mind but, nevertheless, have their own interest. First we approach the 
symmetric cases such as, for example, exponential growths. Second, we approach the 
{\it limit} cases, that is, the ones whose best approximation is not an exponential 
but the limit as $k\to 0$ or $k\to-\infty$ of exponentials. They are not what 
one can expect to find while adjusting data that follow an exponential 
decay but we have been able to identify when they occur and deal 
successfully with them.

\subsection{Symmetric cases}\label{sbs_sym}
If $(\t,T)$ and $f$ are as in the statement of Lemma~\ref{44}, then a moment's 
reflection is enough to realize that:
\begin{enumerate}
\item The $t$-{\em symmetric} data $((-t_4,-t_3,-t_2,-t_1),(T_4,T_3,T_2,T_1))$ 
have 
\begin{equation*}
g_1(t)=f(-t)=a\exp(-kt)+b
\end{equation*}
as its best approximation. 
\item The $T$-{\em symmetric} data $((t_1,t_2,t_3,t_4),(-T_4,-T_3,-T_2,-T_1))$ 
have \begin{equation*}
g_2(t)=-f(t)=-a\exp(kt)-b
\end{equation*} 
as its best approximation. 
\item The {\em bisymmetric} data $((-t_4,-t_3,-t_2,-t_1),(-T_4,-T_3,-T_2,-T_1))$
have \begin{equation*}
g_3(t)=-f(-t)=-a\exp(-kt)-b
\end{equation*}
as its best approximation.
\end{enumerate}
These symmetries correspond to the following:
\begin{enumerate}
\item If $(T_2-T_4)/(t_4-t_2)<(T_1-T_3)/(t_3-t_1)<0$, then there are still two 
changes of sign in the coefficients of $p$, so it has another positive root. 
The difference here is that both $a$ and $k$ are positive. Please observe that 
this means that $(T_2-T_4)/(t_2-t_4)>(T_1-T_3)/(t_1-t_3)>0$, so $f$ must be 
increasing and increases faster for greater $t$. 
\item If $(T_1-T_3)/(t_3-t_1)<(T_2-T_4)/(t_4-t_2)<0$, then $a<0$ and $k<0$. 
\item If $(T_2-T_4)/(t_4-t_2)>(T_1-T_3)/(t_3-t_1)>0$, then everything goes indisturbed but we have 
$p'(1)<0$, so the second root of $p$ is greater than 1. This implies that $k=\exp(z)\in(0,\infty)$ 
and $a<0$. 
\end{enumerate}

\subsection{Limit cases}\label{sbs_limit}
Even if the conditions are not fulfiled by any symmetric version of the dataset, 
the computations made in the proof of Lemma~\ref{44} give the answer to the 
approximation problem: 
\begin{enumerate}
\item If $(T_1-T_3)/(t_3-t_1)=(T_2-T_4)/(t_4-t_2)>0$, then $z=1$ is a double root 
--this corresponds to $k=0$-- and the ``exponential" we are looking for is a line 
with negative slope. Namely, its slope is $(T_3-T_1)/(t_3-t_1)$ and this best 
approximation is the line given by Chebyshev's Alternation Theorem. 
\item If $(T_1-T_3)/(t_3-t_1)=(T_2-T_4)/(t_4-t_2)<0$, then this ``exponential" 
is a line with positive slope $(T_3-T_1)/(t_3-t_1)$ -- this is symmetric to the 
previous case. 
\item If $T_1=T_3$ or $T_2=T_4$, then we have, up to symmetries, three cases: \\ 
\underline{i}: If $T_2=T_4=\m(T)$ and $T_1\leq T_3$, then the best approximation is a constant. 
Namely, $f(t)=\frac 12(T_2+ T_3)$. \\
\underline{ii}: If $T_1>T_3\geq T_2=T_4$ then there is no global best approximation, but every 
exponential approximates $(\t,T)$ worse than the limit, with $k\to-\infty$, 
of the best approximations. This limit is 
\begin{equation}\label{fcte}
\f_{-\infty}(\t)=\lim_{k\to-\infty}\f_k(\t)=
\left(T_1-\frac{T_3-T_2}2,\frac{T_3+T_2}2,\frac{T_3+T_2}2,\frac{T_3+T_2}2\right),
\end{equation}
and it turns out to be also a kind of best approximation for every 
$T_4\in[T_2,T_3]$. \\
\underline{iii}: If $T_1> T_2=T_4>T_3$ then the situation is as follows: \\
As $\M(T)$ lies before $\m(T)$, any good approximation must be nonincreasing. 
$T$ attains its second greatest value after $\m(T)$, so every decreasing 
function approximates $(\t,T)$ worse than the function $\f(\t)$ 
defined as in~\eqref{fcte}. Actually, $(t_2,T_2)$ could be ignored 
whenever $T_2\in[T_3,T_4]$, as we are about to see in the last item: 
\item Finally, if $T_1-T_3$ and $T_2-T_4$ have different signs, then there is just one 
change of signs in the coefficients of $p$, so the only positive root of $p$ 
is $z=1$ and $k=0$, and there is no function fulfilling the statement. More 
precisely, this situation has two paradigmatic examples with 
 $\t=(1,2,3,4)$ and $T=(3,0,1,2)$ or with $\t'=(1,2,3,4)$ and $T'=(3,1,0,2)$. \\
In the first case, the third point $(3,1)$ simply does not affect the 
approximation in the sense that, for every $k$, the exponential $f_k:\R\to\R$ 
that best approximates $T$ fulfils 
\begin{equation*}
\|(f_k(1)-T_1,f_k(2)-T_2,f_k(3)-T_3,f_k(4)-T_4)\|_\infty>|f_k(3)-T_3|.
\end{equation*} 
Namely, if $f_k$ is decreasing then $f_k(4)-T_4<f_k(3)-T_3<f_k(2)-T_2$, so 
$f_k(3)-T_3$ is neither $\M(\f_k(\t)-T)$ nor $\m(\f_k(\t)-T)$. If $f_k$ is increasing then 
$f_k(1)-T_1<f_k(i)-T_i$ for $i=2, 3, 4$ and this, along with Proposition~\ref{rrr}, 
implies that $f_k$ cannot be the best approximation. \\ 
The second case is similar. Though the point $(t_2',T_2')=(2,1)$ is relevant for some 
approximations, it is skippable for every $k<k_0$ for some $k_0<0$. 
\end{enumerate}

\section{General datasets} \label{s_general}
In this Section, we apply the previous results to datasets with arbitrary size 
in order to find out when a dataset has an exponential as its best approximation. 
Before we arrive to this first objective main result, Theorem~\ref{t_discrimina_cuando_ma}, 
we will need several minor results. The path that we will follow is, in a 
nutshell, the following: 

Lemma~\ref{mismostres} is a technical Lemma that allows us to show that 
the maps $k\mapsto a_k, k\mapsto b_k, k\mapsto \f_k(\t)$ are continuous, 
see Corollary~\ref{akbkcontinuas}. 

Lemma~\ref{recta} is just Chebyshev's Alternation Theorem, and it suffices 
to determine which datasets are good to be approximated by decreasing, convex, 
functions --like exponential decays. We will call that datasets 
\textit{admissible} from Definition~\ref{defadmi} on.

Then we determine the vectors one obtains by taking limits of exponentials 
with exponents converging to $\pm\infty$ or 0, see Proposition~\ref{limites}. 

With all these preparations, we are ready to translate Lemma~\ref{44} to a 
more general statement, keeping the $n=4$ condition. We give a necessary and 
sufficient condition for any dataset $(\t,T)$ to be approximable by exponential 
decays in terms that are easily generalizable to $n>4$. 
This is Proposition~\ref{propadmi}. 

In Proposition~\ref{propadmi} and \ref{TablaLim} we improve the results in 
Corollary~\ref{akbkcontinuas} to get Remark~\ref{corobueno}, where we 
show that we can handle the best approximations at ease if the variations 
of $k$ are small enough. 

Finally, Proposition~\ref{p_de_n_a_4}, reduces the general problem to 
the $n=4$ case, thus getting Theorem~\ref{t_discrimina_cuando_ma}.

\begin{Lemma}\label{mismostres}
Let $k_0\neq 0$ and $f_{k_0}(t)=a_{k_0}\exp(kt)+b_{k_0}$ be the best approximation 
for ${k_0}$ and suppose that there are exactly three indices $i<j<m$ such that the equalities 
\begin{equation*}
T_i-f_{k_0}(t_i)=f_{k_0}(t_j)-T_j=T_m-f_{k_0}(t_m)=\delta\|T-\f_{k_0}(\t)\|_\infty,\quad \delta=\pm 1
\end{equation*} 
hold. Then, there exists $\e>0$ such that, for every $k\in({k_0}-\e,{k_0}+\e)$ 
the equalities hold with the same indices. Moreover, if $k'>k_0$ is such that the 
indices where the norm is attained are not $i<j<m$, then there exists 
$k''\in[k_0,k']$ for which the norm is attained in at least four indices. 
\end{Lemma}

\begin{proof}
Suppose $\delta=1$, the case $\delta=-1$ is symmetric. 

As $a_{k_0}=\frac{T_i-T_m}{\exp(k_0t_i)-\exp(k_0t_m)}$, see Lemma~\ref{fijok3r}, 
taking $\alpha_k=\frac{T_i-T_m}{\exp(kt_i)-\exp(kt_m)}$ for every $k\neq 0$ 
we have $\alpha_{k_0}=a_{k_0}$. If we take, further
\begin{equation*}
\beta_k=\frac 12(T_i+T_j-\alpha_k(\exp(kt_i)+\exp(kt_j))),
\end{equation*}
then $\beta_{k_0}=b_{k_0}$, so defining 
$\cl{f_k}(t)=\alpha_k\exp(kt)+\beta_k$ we get $\cl{f_{k_0}}=f_{k_0}$ 
and a straightforward computation shows that 
\begin{equation*}
T_i-\cl{f_{k}}(t_i)=\cl{f_{k}}(t_j)-T_j=T_m-\cl{f_{k}}(t_m)
\end{equation*}
for every $k$. Given any $l\in\{1,\ldots,n\}\setminus\{i,j,m\}$, 
our hypotheses give 
\begin{equation*}
T_j-f_{k_0}(t_j)<T_l-f_{k_0}(t_l)<T_i-f_{k_0}(t_i).
\end{equation*}
As the map 
$k\mapsto \cl{f_{k}}(t_l)=\alpha_{k}\exp(kt_l)+\beta_{k}$ 
is continuous for every $l$, 
we obtain that 
\begin{equation}\label{jli}
T_j-\cl{f_k}(t_j)<T_l-\cl{f_k}(t_l)<T_i-\cl{f_k}(t_i)
\end{equation}
holds for every $k$ in a neighbourhood of $k_0$, say $({k_0}-\e_l,{k_0}+\e_l)$. 
Since there are only finitely many indices we may take $\e$ as the minimum of the 
$\e_l$ to see that $\cl{f_{k}}$ is the best approximation for $k\in({k_0}-\e,{k_0}+\e)$, 
and finish the proof of the first part. 

As for the moreover part, it is quite obvious that the expression for $f_k$ 
will be $f_k(t)=a_k\exp(kt)+b_k$ with 
\begin{equation*}
a_k=\frac{T_i-T_m}{\exp(kt_i)-\exp(kt_m)}\ \mathrm{and\ }
b_k=\frac 12(T_i+T_j-a_k(\exp(kt_i)+\exp(kt_j)))
\end{equation*}
if and only if, for every $l\not\in\{i,j,m\}$, one has 
\begin{equation*}
T_j-\frac{T_i-T_m}{\exp(kt_i)-\exp(kt_m)}t_j \leq 
T_l-\frac{T_i-T_m}{\exp(kt_i)-\exp(kt_m)}t_l \leq 
T_i-\frac{T_i-T_m}{\exp(kt_i)-\exp(kt_m)}t_i.
\end{equation*}
Please observe that the symmetric inequalities could hold and it would 
make $f_k$ have the same expression, but this would imply $n=3$.
Anyway, if $k'>k$ is such that there is some $l$ for which 
\begin{equation*}
T_l-\frac{T_i-T_m}{\exp(k't_i)-\exp(k't_m)}t_l\not \in 
\left[T_j-\frac{T_i-T_m}{\exp(k't_i)-\exp(k't_m)}t_j,T_i-\frac{T_i-T_m}{\exp(k't_i)-\exp(k't_m)}t_i\right],
\end{equation*}
then it is clear that there is $\cl{k}\in[k,k']$ such that 
\begin{equation*}
T_l-\frac{T_i-T_m}{\exp(\cl{k}t_i)-\exp(\cl{k}t_m)}t_l=T_j-\frac{T_i-T_m}{\exp(\cl{k}t_i)-\exp(\cl{k}t_m)}t_j\ \mathrm{or\ }
\end{equation*}
\begin{equation*}
T_l-\frac{T_i-T_m}{\exp(\cl{k}t_i)-\exp(\cl{k}t_m)}t_l=T_i-\frac{T_i-T_m}{\exp(\cl{k}t_i)-\exp(\cl{k}t_m)}t_i.
\end{equation*} 
Maybe it is not $k''=\cl{k}$, but taking $k''$ as the smallest real number in $(k,k']$ 
for which there exists such an $l$, we are done. 
\end{proof}

\begin{Corollary}\label{akbkcontinuas}
The maps $k\mapsto a_k, k\mapsto b_k $ and $k\mapsto \f_k(\t)$ are continuous. 
\end{Corollary}

\begin{Lemma}\label{recta}
Let $T,\t\in\R^n$. There exists exactly one line $r(t)=at+b$ 
such that 
\begin{equation}\label{eqnrecta}
T_i-(at_i+b)=at_j+b-T_j=T_m-(at_m+b)=\delta\|T-(a\t+b)\|_\infty 
\end{equation}
for some $1\leq i<j<m\leq n$ and $\delta=\pm 1$, and this line 
approximates $T$ better than any other line. 
\end{Lemma}

\begin{proof}
It is a particular case of the Chebyshev's Alternation Theorem, applied to the 
polygonal defined by $(\t,T)$.
\end{proof}

\begin{Remark}\label{remadmi}
Thanks to Lemma~\ref{recta}, we can define which vectors $T$ will be our ``good vectors": 
those for which $a<0$ and the equalities~(\ref{eqnrecta}) hold with $\delta=1$
and not with $\delta=-1$. 
When our data fulfil these conditions, we have some idea of decreasing monotonicity 
and also some kind of convexity, and this is the kind of dataset that we wanted, 
though we will need to add some further conditions. 
Anyway, when dealing with datasets that fulfil any couple of symmetric conditions 
we just need to have in mind the symmetries. Namely, they will behave as in 
Subsection~\ref{sbs_sym}. 

% The left-right symmetry corresponds to the variable change $t\mapsto -t$ and the 
% up-down symmetry corresponds to $f\mapsto -f$. 
\end{Remark}

\begin{Definition}\label{defadmi}
Let $r(t)=at+b$ be the line that best approximates $(\t,T)$. 
We will say that $(\t,T)$ is {\em admissible} when $a<0$
and there exist $1\leq i<j<m\leq n$ such that 
\begin{enumerate}
\item 
$T_i-r(t_i)=r(t_j)-T_j=T_m-r(t_m)=\|T-\r(\t)\|_\infty.$
\item $-\|T-\r(\t)\|_\infty\leq r(t_l)-T_l<\|T-\r(\t)\|_\infty$ for every $l<i$ and every $l>m$. 
\item $-\|T-\r(\t)\|_\infty\leq T_l-r(t_l)<\|T-\r(\t)\|_\infty$ for every $i<l<m$. 
\end{enumerate}
\end{Definition}

Once we have stated the kind of data which we will focus on, say {\em discretely 
decreasing and convex}, now we have to determine when they will be approximable. 
Before that, we will study the behaviour of the limits of best approximations. 

\begin{Lemma}\label{lemalimi}
Let $t_1<t_2<t_3$. For each $k\neq 0$ consider some exponential 
$g_k(t)=a_k\exp(kt)+b_k$ and 
\begin{equation*}
\psi(k)=\frac {g_k(t_1)-g_k(t_2)}{g_k(t_2)-g_k(t_3)}.
\end{equation*} 
Then $\psi(k)$ depends on $k$ but not on $a_k$ or $b_k$ and, moreover: 
\begin{enumerate}
\item When $k\to-\infty, \psi(k)\to\infty.$
\item When $k\to 0$, $\lim(\psi(k))=\frac{t_1-t_2}{t_2-t_3}.$
\item When $k\to\infty, \lim(\psi(k))=0$.
\end{enumerate} 
\end{Lemma}

\begin{proof}
We only need to make some elementary computations to show that 
\begin{equation*}
\frac {g_k(t_1)-g_k(t_2)}{g_k(t_2)-g_k(t_3)}=
\frac{\exp(kt_1)-\exp(kt_2)}{\exp(kt_2)-\exp(kt_3)}.
\end{equation*} 
The computation of the limit at 0 only needs a L'Hôpital's rule application, 
the other ones are even easier once one substitutes $z=\exp(k)$. 
See~\cite[Lemma 2.10]{TAC1}. 
\end{proof}

\begin{Proposition}\label{limites}
Let $T,\t\in\R^n$. Then, the following hold:
\begin{enumerate}
\item For $k_0\neq 0$, $\displaystyle\lim_{k\to k_0}\f_k(\t)=\f_{k_0}(\t)$. 
\item For $k_0=0$, $\displaystyle\lim_{k\to k_0}\f_k(\t)$ is $\r(\t)$, the line that 
best approximates $(\t,T)$. 
\item For $k_0=\infty$, $\f_{\infty}(\t)=\displaystyle\lim_{k\to \infty}\f_k(\t)$ takes 
at most two values, and fulfils 
\begin{equation*}
(\f_{\infty}(\t))_1=(\f_{\infty}(\t))_2=\ldots=(\f_{\infty}(\t))_{n-1}.
\end{equation*}
\item For $k_0=-\infty$, $\f_{-\infty}(\t)=\displaystyle\lim_{k\to -\infty}\f_k(\t)$ takes
at most two values, and fulfils 
\begin{equation*}
(\f_{-\infty}(\t))_2=(\f_{-\infty}(\t))_3=\ldots=(\f_{-\infty}(\t))_{n}.
\end{equation*}
\end{enumerate}
\end{Proposition}

\begin{proof}
Let $k_0\neq 0$. If the best approximation for $k_0$ is a constant, then 
it is constant for every $k\neq 0$, and the constants are obviously the same. 
So, we may suppose $f_k$ is not a constant for any $k\neq 0$. 
In this case, Lemma~\ref{mismostres} implies that $k\to f_k(t_l)$ is continuous for 
every $l$. As we have just a finite amount of indices, this means that 
$\displaystyle\lim_{k\to k_0}\f_k(\t)=\f_{k_0}(\t),$ so we are done. 
\footnote{¿No sería conveniente tirar aquí del Corolario 4.2?}

The proof of the three last items is immediate from Lemma~\ref{lemalimi}:
\end{proof}

\begin{Proposition}\label{propadmi}
Let $\t=(t_1,t_2,t_3,t_4)$ and $T=(T_1,T_2,T_3,T_4)$. Then, the best exponential 
approximation to $(\t,T)$ has the form $f(t)=a\exp(kt)+b$ with $a>0, k<0$ if 
and only if $(\t,T)$ is admissible and the following does not happen: 

$\displaystyle\spadesuit$ $T_1=\M(T)$ and the second greatest value 
of $T$ is attained after $\m(T)$. 
\end{Proposition}

\begin{proof}
As the second greatest value of $T$ will appear frequently in this proof, 
we will denote it as $M_T=\max\{T_i:i=2,3,4\}$. Analogously, 
$m_T=\min\{T_i:i=2,3,4\}.$

If $\spadesuit$ happens, then the following is the limit of 
best approximations when $k\to -\infty$
\begin{equation}\label{f-infty}
\f_{-\infty}(\t)=\left(T_1-\frac{(M_T-\m(T))}{2},\frac{(\m(T)+M_T)}{2},
\frac{(\m(T)+M_T)}{2},\frac{(\m(T)+M_T)}{2}\right). 
\end{equation}
Indeed, as $\f_{-\infty}(\t)$ is the pointwise limit of functions fulfiling 
(\ref{eqrrr}), it must fulfil (\ref{eqrrr}) too. It is clear that this implies 
that $\f_{-\infty}(\t)$ must be as in (\ref{f-infty}). 
It is clear that every strictly decreasing function approximates $(\t,T)$ 
worse than $\f_{-\infty}(\t)$, so we have finished the first part of the proof. 

Conversely, if $(\t,T)$ is admissible then there are exponentials 
with $a>0, k<0$ that approximate $(\t,T)$ better than the line $\r(\t)=\f_{0}(\t)$. 
Indeed, we only have to consider the three points of the Definition~\ref{defadmi} 
and take into account Lemma~\ref{mismostres}. 
As the function error is quasiconvex, the only option for contradicting 
the statement is that every exponential is worse than the $-\infty$ limit 
of the approximations $\f_k$, and of course this limit is not better than 
$\f_{-\infty}(\t)$ as in~\eqref{f-infty} because no vector of the form 
$(x,y,y,y)$ approximates $T$ better than this. 
So, we may suppose $\f_{-\infty}(\t)$ is the best approximation --and please 
recall that we are supposing that $(\t,T)$ is admissible. 
We need to break down several possibilities: 

\noindent\ \underline{I}: If $T_1\in[m_T,M_T]$, then we can change the first 
coordinate of $f_{-\infty}(\t)$ from $T_1+(M_T-\m(T))/2$ to $(M_T+\m(T))/2$ 
without increasing the error, so one best approximation is a constant 
and this means that $a=0$, so $(\t,T)$ is not admissible, a contradiction. 

\noindent\ \underline{II}: If $T_1<m_T$ then $(\t,T)$ is not admissible. 

\noindent\ \underline{III}: If $T_1>M_T$, then we still have some options: 

\underline{i}: If $T_2\leq T_4$ then we obtain that $\spadesuit$ holds, no matter 
the value of $T_3$. 

\underline{ii}: If $T_2> T_4$ and the rate of decreasing $(T_2-T_4)/(t_4-t_2)$ 
is greater than $(T_1-T_3)/(t_3-t_1)$, then $(-\t,-T)$ fulfils the hypotheses 
of Lemma~\ref{44}. This implies that the best approximation to $(-\t,-T)$ has 
$a>0, k<0$, so the best approximation to $(\t,T)$ has $a<0, k>0$. 

\underline{iii}: If $T_2> T_4$ and the rates of decreasing are equal, 
then $(\t,T)$ is not admissible because this implies 
\begin{equation*}
T_1-r(t_1)=r(t_2)-T_2=T_3-r(t_3)=r(t_4)-T_4.
\end{equation*}

\underline{iv}: If $T_2> T_4$ and $(T_2-T_4)/(t_4-t_2)<(T_1-T_3)/(t_3-t_1)$ 
then Lemma~\ref{44} ensures that the best approximation is 
$f_k(t)=a\exp(kt)+b$, with $a>0$ and $k<0$. 
\end{proof}

\begin{Notation}\label{TablaLim}
Let $r(t)=T_1-\frac{T_1-T_3}{t_1-t_3}(t-t_1)$ be the line that contains 
$(t_1,T_1)$ and $(t_3,T_3)$. In~\cite{TAC1}, Lemma~2.10, it is seen that, if 
$g_k(t)=a_k\exp(kt)+b$, where 
\begin{equation*}
a_k=\frac{T_1-T_3}{\exp(kt_1)-\exp(kt_3)},\quad \mathrm{and}\quad b=T_1-a_k\exp(kt_1),
\end{equation*} 
then $g_k(t_1)=T_1, g_k(t_3)=T_3$ and 
\begin{itemize}
\item As $k\to\infty$, $g_k(t_4)\to-\infty$. 
\item As $k\to-\infty$, $g_k(t_4)\to T_3$. 
\item As $k\to 0$, $g_k(t_4)\to r(t_4)$ . 
\end{itemize}
Essentially, the same proof suffices to show how $g_k(t_2)$ behaves: 
\begin{itemize}
\item As $k\to\infty$, $g_k(t_2)\to T_1$.
\item As $k\to-\infty$, $g_k(t_2)\to T_3$.
\item As $k\to 0$, $g_k(t_2)\to r(t_2)$ .
\end{itemize}

This implies that the map $k\mapsto g_k(t_2)$ is strictly increasing, while 
$k\mapsto g_k(t_4)$ is strictly decreasing. So, $k\mapsto g_k(t_2)$ increases as 
$k\mapsto g_k(t_4)$ decreases and, moreover, the map $\phi:(T_3,T_1)\to(-\infty,T_3)$ 
given by $\phi(g_k(t_2))= g_k(t_4)$ for $k\neq 0,$ and $\phi(r(t_2))=r(t_4)$ is a 
(decreasing) homeomorphism from $(T_3,T_1)$ to $(-\infty,T_3)$. 
Applying the same reasoning to $t_0<t_1$ and to $t_3<t_4$ we obtain this key result: 
\end{Notation}

\begin{Proposition}\label{decinchom}
Let $t_0<t_1<t_2<t_3<t_4$ and $T_1> T_3$ and consider for every $k\neq 0$ the only 
exponential $g_k(t)=a\exp(kt)+b$ such that $g_k(t_1)=T_1, g_k(t_3)=T_3$ and 
$f_0(t)$ the only line such that $g_0(t_1)=T_1, g_0(t_3)=T_3$. 
Then, all the following maps are homeomorphisms, $\phi_0$ and $\phi_4$ are 
decreasing and $\phi_2$ is increasing: 
\begin{enumerate}
\item $\phi_0:(-\infty,\infty)\to(T_1,\infty)$ defined as $k\mapsto g_k(t_0)$. 
\item $\phi_2:(-\infty,\infty)\to(T_3,T_1)$ defined as $k\mapsto g_{k}(t_2)$. 
\item $\phi_4:(-\infty,\infty)\to(-\infty,T_3)$ defined as $k\mapsto g_k(t_4)$. 
\end{enumerate}
\end{Proposition}

We can rewrite Proposition~\ref{decinchom} as follows: 
\begin{Remark}\label{corobueno}
Let $\t=(t_0,t_1,t_2,t_3,t_4)$ and $\alpha_1>\alpha_3\in\R$. Let, for every 
$k\in\R$, $g_k(t)=c_k\exp(kt)+d_k$ be the only exponential that fulfils 
$g_k(t_1)=\alpha_1, g_k(t_3)=\alpha_3$. 
Then, when $k$ increases, $g_k(t_0)$ and $g_k(t_4)$ decrease and $g_k(t_2)$ 
increases and everything is continuous. 
\end{Remark}

\begin{Proposition} \label{p_de_n_a_4}
The best exponential approximation (including limits) to $(\t,T)$ 
is the best approximation for some quartet 
$((t_{i_1},t_{i_2},t_{i_3},t_{i_4}),(T_{i_1},T_{i_2},T_{i_3},T_{i_4})).$
\end{Proposition}

\begin{proof}
Let $f_{k_0}$ be the best approximation and suppose that the conclusion does not 
hold. Then, we may suppose that there are exactly three indices where the 
norm is attained, say $i<j<m$ and 
\begin{equation*}
T_i-f_{k_0}(t_i)=f_{k_0}(t_j)-T_j=
T_m-f_{k_0}(t_m)=\|T-\f_{k_0}(\t)\|_\infty.
\end{equation*} 
If $\f_{k_0}(\t)$ is the limit at $-\infty$ of the best approximations, 
then it is the best approximation for every quartet that contains 
$((t_i,t_j,t_m),(T_i,T_j,T_m))$ because this means that $\spadesuit$ holds. 
So, suppose that the best approximation 
is $f_{k_0}$, for some $k_0\in\R$ --maybe $k_0=0$. 
Then, for some $\e>0$ the functions $f_k$, with $k\in(k_0-\e,k_0)$ approximate this 
triple better than $f_{k_0}$. Reducing if necessary $\e$, Remark~\ref{corobueno} 
implies that every $f_k$ with $k\in(-\e,0)$ approximates $(\t,T)$ 
better than $f_{k_0}$, thus getting a contradiction. 
\end{proof}

\begin{Theorem}\label{t_discrimina_cuando_ma}
Let $(\t,T)$ be admissible. Then, the best approximation is a exponential if 
and only if $\spadesuit$ does not happen.
\end{Theorem}

\begin{proof}
The proof of Proposition~\ref{propadmi} is enough to see that $\spadesuit$ 
avoids the option of $(\t,T)$ being approximable by a best exponential. 

If $(\t,T)$ is admissible, then the best approximation cannot be the 0-limit 
of exponentials, so it is either an exponential or the $-\infty$-limit of 
exponentials. So, suppose it is the $-\infty$-limit and let us see that in this case 
$\spadesuit$ holds. It is clear that $T_1>M_T=\max\{T_2,\ldots,T_n\}$ as in the proof 
of Proposition~\ref{propadmi}, so just need to show that $M_T$ 
occurs later than $\m(T)$. Let $m_T=\min\{T_2,\ldots,T_n\}$. A moment's 
reflection suffices to realize that $\f_{-\infty}(T_m)=(M_T+m_T)/2$ for every 
$m\geq 2$ and so, the error for $\f_{-\infty}$ is exactly $r=(M_T-m_T)/2$. 
Let $T_i$ be the last appearance of $M_T$ and $T_j$ the first 
appearance of $m_T$ and suppose $i<j$, i.e., that $\spadesuit$ 
does not hold. 
Thanks to Proposition~\ref{decinchom}, for small $\e>0$ there is $k$ close 
enough to $-\infty$ so that we can find $g_k(t)=a\exp(kt)+b$ such that 
$|g_k(t_1)-T_1|<r$, $\f_{-\infty}(T_m)<g_k(t_m)<\f_{-\infty}(T_m)+\e$ when 
$m\in\{2,\ldots,j-1\}$ and $\f_{-\infty}(T_m)>g_k(t_m)>\f_{-\infty}(T_m)-\e$ when 
$m\in\{j,\ldots,n\}$. If we take $\e$ small enough, $g_k$ approximates $(\t,T)$ 
better than $\f_{-\infty}$. 
\end{proof}

% \section{Concluding the computations}\label{s_conclusion}
The value of $b$ in \eqref{ak3bk3} can be easily generalised, so we do not need 
to worry about it. If we are able to determine $k$ and $a_k$, then finding 
$b$ is just a straightforward computation. Namely: 

\begin{Lemma}\label{easyb}
For $k\in(-\infty,0), a\in(0,\infty)$, the best approximation to $T$ in 
$\{a\exp(k\t)+b:b\in\R\}$ is attained when 
\begin{equation}\label{bk}
b=\frac 12(\M(T-a\exp(k\t))+\m(T-a\exp(k\t)).
\end{equation}
With this $b$, the error is 
\begin{equation*}
\frac 12(\M(T-a\exp(k\t))-\m(T-a\exp(k\t)).
\end{equation*}
\end{Lemma}

\begin{proof}
First, we are going to compute the error. Since, obviously, 
\begin{equation*}
\M(|T-a\exp(k\t)-b|)=\max\{\M(T-a\exp(k\t)-b),- \m(T-a\exp(k\t)-b)\},
\end{equation*} 
we just need to take into account that~\eqref{bk} implies 
\begin{equation*}
\begin{split}
&\M(T-a\exp(k\t)-b)=
\M(T-a\exp(k\t))-b= \\ &
\M(T-a\exp(k\t))-\frac 12(\M(T-a\exp(k\t))+\m(T-a\exp(k\t))) = \\
& \frac 12(\M(T-a\exp(k\t))-\m(T-a\exp(k\t))).
\end{split}
\end{equation*}

\begin{equation*}
\begin{split}
&- \m(T-a\exp(k\t)-b)=- \m(T-a\exp(k\t))+b= \\ &
- \m(T-a\exp(k\t))+\frac 12(\M(T-a\exp(k\t))+\m(T-a\exp(k\t)))= \\
& \frac 12(\M(T-a\exp(k\t))-\m(T-a\exp(k\t))).
\end{split}
\end{equation*}
On the one hand, this implies that the error is as in the statement. 
On the other hand, let $b'<b$. Then, 
\begin{equation*}
\M(T-a\exp(k\t)-b')>\M(T-a\exp(k\t)-b),
\end{equation*} 
so $a\exp(k\t)+b'$ approximates $T$ worse than $a\exp(k\t)+b$. 
The same happens with $\m(T-a\exp(k\t)-b')$ if we take 
$b'>b$, so the best approximation is the one with $b$ as in~\eqref{bk}. 
\end{proof}

With this section we have covered the theoretical aspects about the first objective of this paper. Examples in Section \ref{s_examples} are  about Newton Law of Cooling and directly apply what have been developed here, ending this way our first objective.

\section{Examples}\label{s_examples}
In this section we present two different examples. We intend to apply what have been developed in previous sections to fit data following an \textit{exponential} function.

The calculations in this section were carried out by means of a GNU Octave using an AMD Ryzen 7 3700U processor with 16GB of RAM. The system used is an elementary OS 5.1.7 Hera (64-bit) based on Ubuntu 18.04.4 LTS with a Linux kernel 5.4.0-65-generic.

\subsection{Exponential decay in a Newton’s law of cooling process}\label{ss_example_newton}
In section 4 of~\cite{TAC1}, the paper that motivated this one, we presented an example that was the beginning of our work. In this new approach we consider necessary to fit that very same data but using a new tool: approximation through the $\max$-norm.

We are going to fit data coming from a thermometer achieving thermal balance at the bottom of the ocean. The time evolution of the temperature, according to Newton’s law of cooling process, follows an exponential function as
\begin{equation} \label{eqcooling}
P(t)=\lambda_1 e^{kt}+ \lambda_2,
\end{equation}
where $\lambda_1$ and $\lambda_2 \in \mathbb{R}$ and, in the considered case, $k<0$.

We implemented an algorithm to fit, by a pattern as \eqref{eqcooling}, and taking the $\max$-norm as the approximation criteria, the records obtained by the device. The results corresponding to this implementation are gathered in Table \ref{tb_fit_newton}. Table 2 of~\cite{TAC1} shows similar information for the same fit but while using the Euclidean norm, so interested reader can compare both approximations.

\begin{table}[ht]
\begin{center}
\resizebox{0.5\textwidth}{!}{
\begin{tabular}{ccccc}
\toprule 
CPU Time & $k$ & $\lambda_1$ & $\lambda_2$  \\
(in seconds) & & & & \\
\midrule
$0.31419$ & $-0.0026042$ & $5.7259032$ & $-1.3743464$ \\
\bottomrule 
\end{tabular}
}
\caption{Result of implementation of TAC.} 
\label{tb_fit_newton}
\end{center}
\end{table}

\begin{figure}
     \centering
     \begin{subfigure}[b]{0.48\textwidth}
        \centering
        \includegraphics[width=0.7\linewidth]{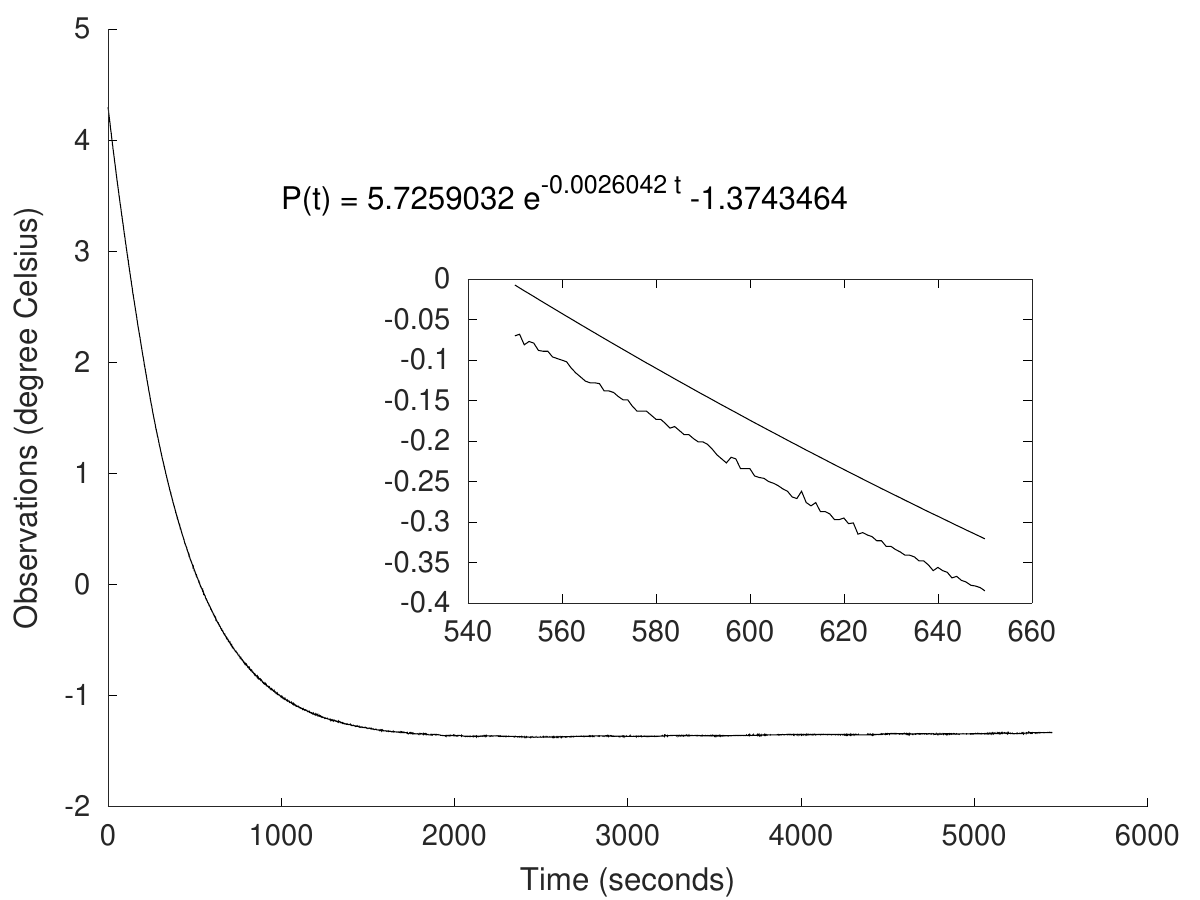}
        \caption{Data and a small detail of TAC fit.}
        \label{subfig:ex_termometro_observaciones}
     \end{subfigure}
     \hfill
     \begin{subfigure}[b]{0.48\textwidth}
        \centering
        \includegraphics[width=0.7\linewidth]{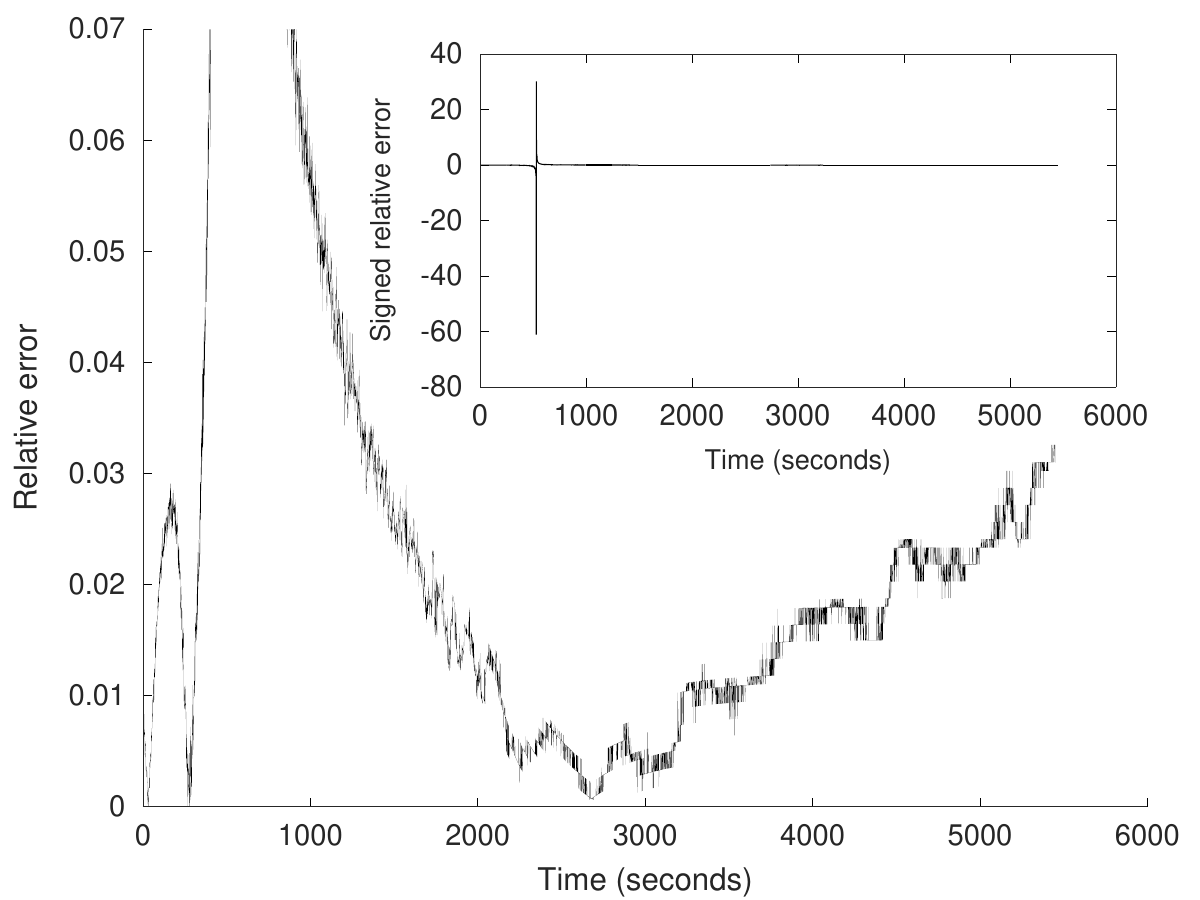}
        \caption{Relative error for values lower or equal than 0.07. on a small scale, the signed relative error.}
        \label{subfig:ex_termometro_relative_error}
     \end{subfigure}
        \caption{Some graphical aspects about this TAC implementation. In the numerical analisys bibliography, relative error is defined with or without sign; in this paper we will consider the latter. A spike can be seen in the small window of Figure \ref{subfig:ex_termometro_relative_error}. This spike should not be considered as an indicator of a poor adjustment of the curve to the data. On the contrary: the spike is due to the proximity of the data to zero and, however, the error remains bounded. This is because curve and data are close enough to control the fact that we are
virtually dividing by zero}
        \label{fig:ex_termometro}
\end{figure}

In Figure \ref{fig:ex_termometro} we present some graphical information about this approximation. Subfigure \ref{subfig:ex_termometro_observaciones} shows observations and fit and Subfigure \ref{subfig:ex_termometro_relative_error} shows the relative error in this implementation. We have designed this figures to resemble figure 3 of~\cite{TAC1} so the reader can graphically compare the results of both approximations. Taking into account each example use a different norm, differences in the fit must exist. Nevertheless they both give a more than reasonable approximation.

\subsection{Exponential decay attaining the max-norm in 4 indices}\label{ss_example_4R}
\begin{figure}[ht]
  \centering
  \includegraphics[width=0.7\linewidth]{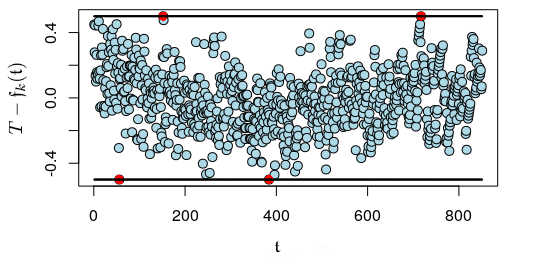}
\caption{Blue dots represent differences between $T$ and $\f_k(\t)$. In 4 dots, $\M(T-\f_k(\t))$ and $\m(T-\f_k(\t))$ are reached, and we have colored them in red. Those are the 4 points mentioned in Proposition \ref{p_de_n_a_4}. Please observe how the maxima and the minima are alternatively reached.}
\label{fig:ex_residuos_4R}
\end{figure}

We wanted to show an example where Proposition \ref{p_de_n_a_4} was fulfilled. Data from Example \ref{ss_example_newton} was not good enough for us to find the 4 points at witch the norm attain its maximum, so we have chosen a different set of data -also coming from a thermometer achieving thermal balance. This new set of data allow us to show a visual perspective of Proposition \ref{p_de_n_a_4}: Figure \ref{fig:ex_residuos_4R} presents in red where the 4 points are.

\section{Economical models and \textit{nlstac}}\label{s_examples_economy}
Now we will cover the second objective: to exemplify the interest for the economic science of the algorithm presented in~\cite{TAC1} and implemented in the R-package \textit{nlstac}. The calculations in this section were carried out by the same machine as in Section \ref{s_examples} using RStudio instead of GNU Octave.

\subsection{The exponential model in demand curves}\label{ss_demand_curves}
In this example we will use the exponential model to fit demand curves. As stated in~\cite{example_koffarnus_a_modified_exponential}, 'behavioral economic demand analyses describe the relationship between the price (including monetary cost and/or effort) of a commodity and the amount of that commodity that is consumed. Such analyses have been successful in quantifying the reinforcing efficacy of commodities including drugs of abuse, and have been shown to be related to other markers of addiction.'

Different mathematical representations of demand curves have been proposed --see, for example,~\cite{example_hursh_economic_demand}. The most widely used model for demand curves in addiction research is presented in the following equation

\begin{equation}\label{eq_demand}
 \log_{10}Q= \log_{10}(Q_0)+k(e^{-\alpha\cdot Q_0 \cdot C}-1),
\end{equation}
where $Q$ represents consumption at a given price, $Q_0$ is known as derived demand intensity, $k$ is a constant that denotes the range of consumption values in log units, $C$ is the commodity price and $\alpha$ is the derived essential value, a measure of demand elasticity.

This model was established by Hursh and Silberberg in~\cite{example_hursh_economic_demand} and has been used in many other works such as~\cite{example_koffarnus_a_modified_exponential}, \cite{example_christensen_demand_for_cocaine} or~\cite{example_strickland_comparing_exponential}. The parameters to be estimated are $Q_0$, $k$ and $\alpha$.

Renaming $\log_{10}(Q_0) -k$ as $b$, $k$ as $a$ and $-\alpha\cdot Q_0$ as $d$, the pattern is now similar to the one we have been working in this paper:
\begin{equation}\label{eq_demand_addapted_renamed}
 \log_{10}Q= a \cdot e^{d \cdot C} + b
\end{equation}

Therefore we simply need to adjust consumption data to the pattern presented in \eqref{eq_demand_addapted_renamed} and undo the changes, being $k=a$, $Q_0=10^{b+k}$ and $\alpha=-\frac{d}{Q_0}$ the parameters we originally sought.

We ran a simulation for this kind of data and successfully fitted it using R package \textit{nlstac}. For the simulation we have established 15 as the number of observations, $Q_0$ as 48, $\alpha$ as 0.006, $k$ as 3.42 and added some noise to the pattern with mean $0$ and standard deviation $0.1$. Tolerance value was set as $10^{-7}$. %Our code is presented in 

Table \ref{tb_fit_consumption} shows the result of this TAC implementation. As can be seen, output is reasonably similar to the values we originally established. Please take into account we added some noise to the data so differences in the parameters were expected.

\begin{table}[ht]
\begin{center}
\resizebox{0.7\textwidth}{!}{
\begin{tabular}{ccccccc}
\toprule 
CPU Time & $Q_0$ & $k$ & $\alpha$ & RSS & MSE  \\
 (in seconds) & & & & \\
\midrule
$0.345$ & $47.74341$ & $3.479838$ & $0.005514662$ & $0.0799227$ & $0.00532818$ \\
\bottomrule 
\end{tabular}
}
\caption{Result of implementation of \textit{nlstac} for pattern \eqref{eq_demand}. RSS denote residual sum of squares and $MSE=RSS/n$ where $n=15$ is the number of observations.} 
\label{tb_fit_consumption}
\end{center}
\end{table}

\begin{figure}[ht]
  \centering
  \includegraphics[width=0.7\linewidth]{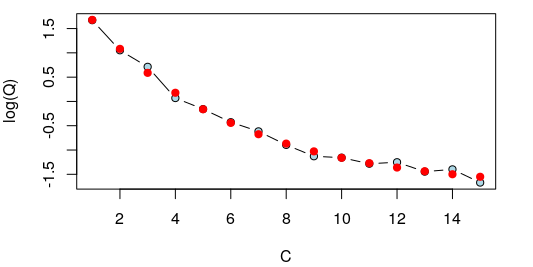}
\caption{Observations in blue, approximation in red.}
\label{fig:ex_consumption}
\end{figure}

In Figure \ref{fig:ex_consumption} we can see the observations (blue dots) and the approximation (red dots). This approximation makes sense: it gets in the middle of the observations, keeping the errors under control.

The calculations in this section were carried out by the same machine as in Section \ref{s_examples} using RStudio instead of GNU Octave.

\subsection{The exponential autoregressive model}\label{ss_ExpAR}
As stated in~\cite{example_xu_modeling_a_nonlinear}, 'nonlinear time series models can reveal nonlinear features of many practical processes, and they are widely used in finance, ecology and some other fields.' Out of those nonlinear time series models, the exponential autoregressive (ExpAR) model is specially relevant. Given a time series $\{ x_1,x_2, x_3,\ldots \}$, the ExpAR model is defined as
\begin{equation*}\label{eq_ExpAR}
 x_t= \left[ \sum_{i=1}^p \left(c_i + \pi_i e^{-\gamma x_{t-1}^2}\right)x_{t-i} \right] + \varepsilon_t,
\end{equation*}
where $\varepsilon_t$ is an i.i.d random variable and independent with $x_i$ and $p$ denotes the system degree and $c_i$, $\pi_i$ (for $i=1,\ldots,p$) and $\gamma$ are the parameters to be estimated from observations. This model can be found in, for example,~\cite{example_xu_modeling_a_nonlinear} or~\cite{example_chen_generalized_exponential_autoregressive}. We followed the notation of the former.

Some generalizations for this model have been made and~\cite{example_chen_generalized_exponential_autoregressive} presents a wide variety of those generalizations. Teräsvirta’s model is an extension of the ExpAR model presented in~\cite{example_Terasvirta} and used in~\cite{example_chen_generalized_exponential_autoregressive}. We will focus in a generalization of Teräsvirta’s model that can be found in equation 10 of~\cite{example_chen_generalized_exponential_autoregressive}:
\begin{equation}\label{eq_ExpAR_Terasvirta_generalization}
 x_t= c_0 + \left[ \sum_{i=1}^p \left(c_i + \pi_i e^{-\gamma (x_{t-d}-z_i)^2}\right)x_{t-i} \right] + \varepsilon_t,
\end{equation}
where $z_i$ (for $i=1,\ldots,p$) are scalar parameters and $d$ is an integer number.

We intend to fit data following \eqref{eq_ExpAR_Terasvirta_generalization} in the particular case when $p=d=2$.  Please observe that this problem is way beyond the proven convergence of TAC algorithm. It requires more than just fitting  of a curve following some \textit{exponential} function since now we have no function to be fitted because every observation depends on the previous ones. This obstacle can be overcome by looking at the problem not as a one-dimensional problem but as a two-dimensional one: if $(x_1,\ldots,x_n)$ are the observations, denoting $y=(x_3,\ldots,x_n)$, $x^1=(x_2,\ldots,x_{n-1})$, $x^2=(x_1,\ldots,x_{n-2})$, $\mathbf{1}=(1,\stackrel{n}{\ldots},1)$, $z^1=z_1 \cdot \mathbf{1}$ and $z^2=z_2 \cdot \mathbf{1}$, data $y$ will depend on two independent variables, $x^1$ and $x^2$, and could be written as:

%$\alpha^0=\alpha_0 \cdot \mathbf{1}$, $\alpha^1=\alpha_1 \cdot \mathbf{1}$, $\alpha^2=\alpha_2 \cdot \mathbf{1}$, $\beta^1=\beta_1 \cdot \mathbf{1}$ and $\beta^2=\beta_2 \cdot \mathbf{1}$

\begin{equation}\label{eq_ExpAR_bidimensional}
 y = c_0 + c_1x^1 + \pi_1 x^1 \sbullet[.5] e^{-\xi (x^2-z^1)^2} + c_2 x^2 + \pi_2 x^2 \sbullet[.5] e^{-\xi (x^2-z^2)^2},
\end{equation}
where operator $\sbullet[.5]$ represents the product between two vectors coordinate to coordinate, that is, given $(a_1,\ldots,a_n),(b_1,\ldots,b_n)\in\mathbb{R}^n$, $(a_1,\ldots,a_n) \sbullet[.5] (b_1,\ldots,b_n):=(a_1 \cdot b_1,\ldots,a_n \cdot b_n)$, and where the exponential and power functions are applied to each coordinate.

As indicated before, this new approach is far for proven in the TAC convergence, however running \textit{nlstac} package will provide us a result that stands to reason. 

For the simulation we have considered in \eqref{eq_ExpAR_Terasvirta_generalization} $p=d=2$. We have generated the $100$ first elements of the time series setting $c_0$ as -1.49, $c_1$ as 1.65, $c_2$ as 0.54, $\pi_1$ as -0.44, $\pi_2$ as -0.84, $\gamma$ as 1.3, $z_1$ as 2.52, $z_2$ as 3.86, $x_1$ as 2.75, $x_2$ as 3.1, and tolerance as $10^{-7}$. %Code is presented in Appendix \ref{anexo_ExpAR}.

In Tables \ref{tb_fit_ExpAR_1} and \ref{tb_fit_ExpAR_2} we gather the results of this implementation. As can be seen, parameters are quite similar to the ones we have previously established.

\begin{table}[ht]
\begin{center}
\resizebox{0.65\textwidth}{!}{
\begin{tabular}{cccccc}
\toprule 
CPU Time & $c_0$ & $c_1$ & $c_2$ & $\pi_1$ & $\pi_2$ \\
 (in seconds) & & & & \\
\midrule
$27.574$ & $-1.5188868$ & $1.6530616$ & $0.5556431$ & $-0.4445339$ & $-0.8514089$ \\
\bottomrule 
\end{tabular}
}
\caption{Result of implementation of \textit{nlstac} for pattern \eqref{eq_ExpAR_Terasvirta_generalization}.} 
\label{tb_fit_ExpAR_1}
\end{center}
\end{table}

\begin{table}[ht]
\begin{center}
\resizebox{0.65\textwidth}{!}{
\begin{tabular}{ccccc}
\toprule 
 $\gamma$ & $z_1$ & $z_2$ & RSS & MSE  \\
\midrule
$1.2729767$ & $2.5142174$ & $3.8667688$ & $2.157779e-06$ & $2.201815e-08$ \\
\bottomrule 
\end{tabular}
}
\caption{Result of implementation of \textit{nlstac} for pattern \eqref{eq_ExpAR_Terasvirta_generalization}, RSS being residual sum of squares and $MSE=RSS/n$ where $n=98$ is the number of observations.} 
\label{tb_fit_ExpAR_2}
\end{center}
\end{table}

\begin{figure}[ht]
  \centering
  \includegraphics[width=0.7\linewidth]{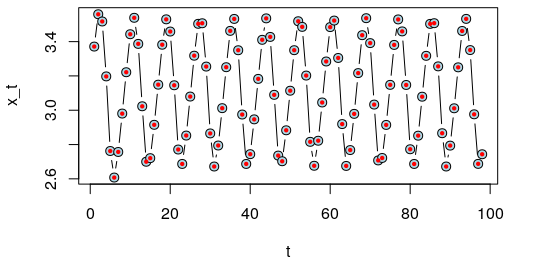}
\caption{Observations in blue, approximation in red (smaller circle). Please observe how each approximation lays over the actual observation.}
\label{fig:ex_ExpAR}
\end{figure}

In Figure \ref{fig:ex_ExpAR} we can see the approximations (red dots) over the actual observations (blue dots, a bit bigger than the red ones), which indicate the approximation is good.

This example is specially relevant because it show us the possibility to use \textit{nlstac} in a problem way different than the one it was intended, so it opens the door to its use in different nonlinear time series models or even in some other models.

%%%%%%%%%%%%%%%%%%%%%%%%%%%%%%%%%%%%%%%%%%

%\vspace{6pt} 

%%%%%%%%%%%%%%%%%%%%%%%%%%%%%%%%%%%%%%%%%%
%% optional
%\supplementary{The following are available online at \linksupplementary{s1}, Figure S1: title, Table S1: title, Video S1: title.}


%%%%%%%%%%%%%%%%%%%%%%%%%%%%%%%%%%%%%%%%%%
\authorcontributions{As can be deduced by the acknowledgments, all authors have been involved during the whole creative process of this paper. However we would like to note that the contributions of Javier Cabello Sánchez have been crucial in the results involving the $max$-norm. This same comment can be suitable for the examples presented in this paper, where Juan Antonio Fernández Torvisco has had a predominant role. Nevertheless Mariano R. Arias has been key to obtain the final version of the paper; connecting both major goals of this paper.}

\funding{This research was partially supported by MICINN [research project references CTM2010-09635 (subprogramme ANT) and PID2019-103961GB-C21], MINECO [research project reference MTM2016-76958-C2-1-P], and Consejería de Economía e Infraestructuras de la Junta de Extremadura [research projects references: IB16056 and GR15152]}.

\institutionalreview{Not applicable}

\informedconsent{Not applicable}

\dataavailability{In this section, please provide details regarding where data supporting reported results can be found, including links to publicly archived datasets analyzed or generated during the study. Please refer to suggested Data Availability Statements in section ``MDPI Research Data Policies'' at \url{https://www.mdpi.com/ethics}. You might choose to exclude this statement if the study did not report any data.} 

\acknowledgments{We thank the workmates that always had interesting comments or hints during coffee breaks. Each author is grateful to the rest for their good attitude which have promoted a positive working atmosphere. The endurance of Javier Cabello Sánchez has been key in this work, and the other authors are grateful to him.}

\conflictsofinterest{The authors declare no conflict of interest.} 

%% Optional
%\sampleavailability{Samples of the compounds ... are available from the authors.}

%%%%%%%%%%%%%%%%%%%%%%%%%%%%%%%%%%%%%%%%%%
%% Only for journal Encyclopedia
%\entrylink{The Link to this entry published on the encyclopedia platform.}

%%%%%%%%%%%%%%%%%%%%%%%%%%%%%%%%%%%%%%%%%%
%% Optional
\abbreviations{The following abbreviations are used in this manuscript:\\

\noindent 
\begin{tabular}{@{}ll}
ExpAR & Exponential autoregressive model
\end{tabular}}

%%%%%%%%%%%%%%%%%%%%%%%%%%%%%%%%%%%%%%%%%%
%% Optional
%\appendixtitles{no} % Leave argument "no" if all appendix headings stay EMPTY (then no dot is printed after "Appendix A"). If the appendix sections contain a heading then change the argument to "yes".
%\appendixstart
%\appendix
%\section{}
%\subsection{}
%The appendix is an optional section that can contain details and data supplemental to the main text---for example, explanations of experimental details that would disrupt the flow of the main text but nonetheless remain crucial to understanding and reproducing the research shown; figures of replicates for experiments of which representative data are shown in the main text can be added here if brief, or as Supplementary Data. Mathematical proofs of results not central to the paper can be added as an appendix.
%
%\begin{specialtable}[H] 
%%\tablesize{\scriptsize}
%\caption{This is a table caption. Tables should be placed in the main text near to the first time they are~cited.\label{taba1}}
%%\tablesize{} % You can specify the fontsize here, e.g., \tablesize{\footnotesize}. If commented out \small will be used.
%\begin{tabular}{ccc}
%\toprule
%\textbf{Title 1}	& \textbf{Title 2}	& \textbf{Title 3}\\
%\midrule
%Entry 1		& Data			& Data\\
%Entry 2		& Data			& Data\\
%\bottomrule
%\end{tabular}
%\end{specialtable}

%\section{}
%All appendix sections must be cited in the main text. In the appendices, Figures, Tables, etc. should be labeled, starting with ``A''---e.g., Figure A1, Figure A2, etc. 

%%%%%%%%%%%%%%%%%%%%%%%%%%%%%%%%%%%%%%%%%%
\end{paracol}
\reftitle{References}


\externalbibliography{yes}
\bibliography{bibliografia_quasiconvex}

\end{document}